\documentclass{article}
\usepackage{amssymb}
\newtheorem{definition}{Definition}
\newtheorem{proposition}{Proposition}
\newtheorem{theorem}{Theorem}

\begin{document}

\begin{titlepage}

\title{Quantization of Teichm\"uller spaces and
the quantum dilogarithm}

\author{R.M. Kashaev\thanks{On leave of absence from
St. Petersburg Branch of the Steklov Mathematical Institute}\\ \\
Helsinki Institute of Physics\\
P.O.Box 9, (Siltavuorenpenger 20 C)\\
FIN-00014 University of Helsinki, Finland\\}

\date{May 1997}

\maketitle

\begin{abstract}
The Teichm\"uller space of punctured surfaces with
the Weil--Petersson symplectic structure and the action of the mapping
class group
is realized as the Hamiltonian reduction of a finite dimensional
symplectic space where the mapping class group acts by
symplectic rational transformations. Upon quantization
the corresponding (projective) representation of the mapping class
group is generated by the quantum dilogarithms.
\end{abstract}

\end{titlepage}

\section*{Introduction}
Despite the progress in understanding and solving the quantum
Chern--Simons field theory with a compact gauge group
(Witten 1989) much remains
unclear in the case of non-compact groups. Among physical
motivations for the study
of the latter case are the interpretation of
$2+1$--dimensional gravity as the Chern--Simons 
theory with a non-compact gauge group (Witten 1988/89) and
identification of the Hilbert space of physical states of the 
Chern--Simons theory with $SL(2,\mathbb{R})$ gauge 
group with the space of Virasoro conformal blocks (Verlinde 1990).
The major mathematical motivation is provided by the possibility of
constructing new topological three--manifold invariants
(Witten 1988/89, 1989).

The purpose of this paper is to quantize the Teichm\"uller
space  of punctured surfaces with the Weil--Petersson simplectic structure.
This space can be identified with (a part of) the phase space of the 
$SL(2,\mathbb{R})$ Chern--Simons theory (Verlinde 1990; Witten 1990).
We start from the Penner parameterization
of the (decorated) Teichm\"uller space (Penner 1987), 
where the mapping class group
is realized explicitly through rational transformations generated by
compositions of the elementary Ptolemy transformation. The latter
transformation is canonical and the quantum dilogarithm
(Faddeev and Kashaev 1992; Faddeev 1995)
implements this transformation on the quantum level. 
Our approach is similar
to the combinatorial quantization of the Chern--Simons theory with
compact gauge groups (Fock, Rosly 1992; Alekseev {\em et al.} 1995;
Buffenoir, Roche 1996).

The paper is organized as follows. In Sect.~\ref{section1} first
the results of Penner, to be used in the paper,
are formulated using the
category language, and then the Teichm\"uller space
is described as the Hamiltonian reduction of a 
finite dimensional symplectic space. The quantization
is performed in Sect.~\ref{section2}. All the proofs are omitted
for they are straightforward verifications.

{\bf Acknowledgements}. It is a pleasure to thank N.M. Bogoliubov,
L.D. Faddeev, V.V. Fock, M.A. Semenov-Tian-Shansky,
O. Tirkkonen, A.Yu. Volkov for stimulating discussions.

\section{The phase space}\label{section1}
\subsection{Notation}

In this paper $\mathbb{R}_+$ denotes the set of {\it strictly} positive
real numbers with a Lie group structure given by the multiplication.
If $X$ is a finite set, 
$\mathbb{R}_+^X:=\{f\colon X\to \mathbb{R}_+\}$ denotes the 
space of positive functions on $X$
which is an Abelian group w.r.t. the pointwise multiplication.
The corresponding factor--group w.r.t. the subgroup of constant
functions is denoted
\[
P\mathbb{R}_+^X:=\mathbb{R}_+^X/\mathbb{R}_+.
\]
Each element $x\in X$ is identified
with a coordinate function on $\mathbb{R}_+^X$,
and projective coordinate function in $P\mathbb{R}_+^X$:
\[
X\ni x\colon \mathbb{R}_+^X\ni f\mapsto x(f):=
f(x)\in \mathbb{R}_+.
\]
For set $Y$ groupoid $\mathfrak{G}_Y$ is defined by
$\mathrm{Ob}\mathfrak{G}_Y:=Y$, each set of morphisms
consisting of only one element,
\[
\mathrm{Mor}(x,y):=\{f\colon y\to x\}=\{x\cdot y\},\quad\forall x,y\in Y.
\]
$\mathfrak{P}$ denotes a category with 
$\mathrm{Ob}\mathfrak{P}:=\{P(M,G,\pi,\alpha)\}$, where
$P(M,G,\pi,\alpha)$ is a principal bundle 
with total
space $P$, base manifold $M$, structure group $G$, projection
$\pi\colon P\to M$, and structure group invariant closed two-form 
$\alpha$ in $P$.
The set $\mathrm{Mor}(P_1,P_2)$, where
$P_i:=P_i(M_i,G_i,\pi_i,\alpha_i)$, $i=1,2$, consists of principal bundle
morphisms $f\colon P_2\to P_1$ such that
$f^*\alpha_1=\alpha_2$.

\subsection{The decorated Teichm\"uller space}

Let $\Sigma$
be a closed oriented surface of genus $g$  with  
$s>0$ removed points (punctures) $P_1,\ldots,P_s$, where
\[
2g-2+s>0.
\]
Denote the set of punctures
\[
V:=\{P_1,\ldots,P_s\}.
\]
Let $\mathcal{T}_\Sigma$ be the Teichm\"uller space of
marked conformal types
of hyperbolic metrics on $\Sigma$, and $\tilde\mathcal{T}_\Sigma$, the
decorated Teichm\"uller space, which is
a principal $\mathbb{R}_+^s$ foliated fibration
$\phi\colon\tilde{\mathcal{T}}_\Sigma\to\mathcal{T}_\Sigma$,
where the fiber over a point of $\mathcal{T}_\Sigma$ is the space of all 
horocycles about the punctures of $\Sigma$.
\begin{definition}
A homotopy class of a path, running between $P_i$ and $P_j$, is called
{\rm ideal arc (i.a.)}. 
A set of ideal arcs, obtained by taking family $X$ of disjointly embedded 
simple arcs in $\Sigma$ running between punctures and subject to the condition
that each component of $\Sigma\setminus X$ is a triangle, 
is called {\rm ideal triangulation (i.t.)}. The set of all i.t. 
on $\Sigma$ is denoted $\Delta_\Sigma$.
\end{definition}
Suppose that $a,b,c,d,e\in\tau\in\Delta_\Sigma$, 
are such that $\{a,b,e\}$ 
and $\{c,d,e\}$ bound distinct triangles. 
The operation of changing i.t.
$\tau$ into i.t. $\tau^e$, which consists in replacing 
i.a. $e$
by a complementary i.a. $e'$ such that triangles, bounded by
$\{a,b,e\}$ and $\{c,d,e\}$,
are replaced by triangles, bounded by $\{b,c,e'\}$ and $\{d,a,e'\}$,
is called {\em elementary move} along i.a. $e$, see Fig.~\ref{figure1}. 
\begin{figure}[h]
\centering
\begin{picture}(200,40)
\put(20,0){\line(-1,1){20}}\put(5,3){$c$}
\put(20,0){\line(1,1){20}}\put(32,3){$b$}
\put(0,20){\line(1,1){20}}\put(5,32){$d$}
\put(40,20){\line(-1,1){20}}\put(32,32){$a$}
\put(20,0){\line(0,1){40}}\put(22,18){$e$}
\put(180,0){\line(-1,1){20}}\put(165,3){$c$}
\put(180,0){\line(1,1){20}}\put(192,3){$b$}
\put(160,20){\line(1,1){20}}\put(165,32){$d$}
\put(200,20){\line(-1,1){20}}\put(192,32){$a$}
\put(160,20){\line(1,0){40}}\put(178,22){$e'$}
\put(20,0){\circle*{3}}
\put(0,20){\circle*{3}}
\put(20,40){\circle*{3}}
\put(40,20){\circle*{3}}
\put(180,0){\circle*{3}}
\put(160,20){\circle*{3}}
\put(180,40){\circle*{3}}
\put(200,20){\circle*{3}}
\put(90,17){$\longrightarrow$}
\end{picture}
\caption{}\label{figure1}
\end{figure}

To each $\tau\in\Delta_\Sigma$ associate object
\[
R(\tau):=\mathbb{R}_+^{\tau}(\mathbb{R}_+^{\tau}/\mathbb{R}_+^V,
\mathbb{R}_+^V,\pi_{R(\tau)},\alpha_\tau)\in\mathrm{Ob}\mathfrak{P},
\]
where the free action of the structure group is defined by
\begin{equation}\label{group-action}
\mathbb{R}_+^V\ni f\colon\mathbb{R}_+^{\tau}\to
\mathbb{R}_+^{\tau},\quad f^*(c):=f(c_0)f(c_1)c,\quad c\in\tau,
\end{equation}
with $c_0$ and $c_1$ being the punctures connected by i.a. $c$; 
\begin{equation}\label{two-form}
\alpha_\tau:=\sum d\ln a\wedge d\ln b+
d\ln b\wedge d\ln c+d\ln c\wedge d\ln a,
\end{equation}
where summation is taken over all triangles with $a,b,c$
being edges of a triangle taken in the clockwise order 
w.r.t. the orientation of $\Sigma$. Define also
\[
R(\star):=\tilde{\mathcal{T}}_\Sigma(\mathcal{T}_\Sigma,\mathbb{R}_+^s,
\phi,\phi^*\omega_{WP}/2)\in\mathrm{Ob}\mathfrak{P},
\]
where $\omega_{WP}$ is the Weil--Petersson K\"ahler form in 
$\mathcal{T}_\Sigma$.

For each $\tau\in\Delta_\Sigma$ define mapping
\begin{equation}\label{l-map}
R(\tau\cdot\star)\colon R(\star)\to
R(\tau),\quad
R(\tau\cdot\star)(h)\colon\tau\ni c\mapsto 
\sqrt{2}e^{\delta_h(c)/2}\in\mathbb{R}_+,
\quad h\in\tilde{\mathcal{T}}_\Sigma,
\end{equation}
where real number 
$\delta_h(c)$ is the signed $\phi(h)$-Poincar\'e distance  between 
the horocycles about the puncture(s) connected by $c$ 
along the geodesic isotopic to $c$.

For a pair of i.t $\tau$ and $\tau^e$, connected by the elementary
move along i.a. $e$, associate mapping
\begin{equation}\label{r-map}
R(\tau^e\cdot\tau)\colon R(\tau)\to R(\tau^e),\quad
R(\tau^e\cdot\tau)^*(x):=\left\{\begin{array}{cl}
(ac+bd)/e&\mathrm{if}\ x=e',\\
x&\mathrm{otherwise;}
\end{array}\right.
\end{equation}
see Fig.~\ref{figure1} for the notation of i.a. 
\begin{theorem}[Penner 1987a, 1987b]
$R$ extends to a unique covariant functor from 
$\mathfrak{G}_{\Delta_\Sigma\cup\star}$ into $\mathfrak{P}$.
\end{theorem}
The mapping class group $M_\Sigma$ of 
$\Sigma$ naturally acts both in $\tilde{\mathcal{T}}_\Sigma$ and
$\Delta_\Sigma$, each i.t.
as a set being mapped to its image in $\Delta_\Sigma$. 
To each $m\in M_\Sigma$ and 
$x\in\Delta_\Sigma\cup\star$ associate
morphism $\mathfrak{r}_m(x)\in \mathrm{Mor}(R(x),R(x))$ such that
\[
R(x)\ni f\mapsto\mathfrak{r}_m(x)(f):=\left\{
\begin{array}{cl}
R(x\cdot m(x))(f\circ m^{-1})& x\in\Delta_\Sigma;\\
m(f)&x=\star.
\end{array}\right.
\]
\begin{theorem}[Penner 1987a]
Mapping $m\mapsto\mathfrak{r}_m$ is a group homomorphism
from $M_\Sigma$ into functorial isomorphisms from $R$ to $R$.
\end{theorem}

\subsection{Teichm\"uller space as the phase space of a constraint system}
\begin{definition}
An i.t. with a choice of a distinguished corner for each triangle is
called {\rm decorated ideal triangulation (d.i.t.)}. The set
of all d.i.t. of $\Sigma$ is denoted $\tilde{\Delta}_\Sigma$.
\end{definition}
There is a natural projection functor
\[
E\colon\mathfrak{G}_{\tilde{\Delta}_\Sigma}
\to\mathfrak{G}_{\Delta_\Sigma}.
\]
Denote by $\dot{\tau}$ the corresponding to 
$\tau\in\tilde{\Delta}_\Sigma$ set of
triangles on $\Sigma$ with distinguished corners.

D.i.t. $\tau_t$, obtained from d.i.t $\tau$ 
by a change of the distinguished corner 
of triangle $t$ as indicated in Fig.~\ref{figure2},
\begin{figure}[h]
\centering
\begin{picture}(200,40)
\put(0,0){\begin{picture}(40,20)
\put(0,0){\line(1,0){40}}
\put(0,0){\line(1,1){20}}
\put(20,20){\line(1,-1){20}}
\put(0,0){\circle*{3}}
\put(20,20){\circle*{3}}
\put(40,0){\circle*{3}}
\put(33,0){$*$}
\put(18,5){$t$}
\end{picture}}
\put(160,0){\begin{picture}(40,20)
\put(0,0){\line(1,0){40}}
\put(0,0){\line(1,1){20}}
\put(20,20){\line(1,-1){20}}
\put(0,0){\circle*{3}}
\put(20,20){\circle*{3}}
\put(40,0){\circle*{3}}
\put(17.5,14){$*$}
\put(18,5){$t'$}
\end{picture}}
\put(95,8){$\longrightarrow$}
\end{picture}
\caption{}\label{figure2}
\end{figure}
is said to be obtained from $\tau$ by 
the {\it elementary change of decoration}
in triangle $t$.
D.i.t. $\tau^e$, obtained from d.i.t. $\tau$
by the elementary move
along i.a. $e$, where distingushed corners of the surrounding triangles
are such as indicated in Fig.~\ref{figure3}, 
\begin{figure}[h]
\centering
\begin{picture}(200,40)
\put(0,0){
\begin{picture}(40,40)
\put(20,0){\line(-1,1){20}}
\put(40,20){\line(-1,-1){20}}
\put(0,20){\line(1,1){20}}
\put(40,20){\line(-1,1){20}}
\put(20,0){\line(0,1){40}}
\put(8,18){$x$}\put(26,18){$y$}
\put(20,0){\circle*{3}}
\put(0,20){\circle*{3}}
\put(20,40){\circle*{3}}
\put(40,20){\circle*{3}}
\put(16,33){$*$}
\put(19.5,33){$*$}
\end{picture}}
\put(160,0){
\begin{picture}(40,40)
\put(20,0){\line(-1,1){20}}
\put(40,20){\line(-1,-1){20}}
\put(0,20){\line(1,1){20}}
\put(40,20){\line(-1,1){20}}
\put(40,20){\line(-1,0){40}}
\put(17,26){$x'$}\put(17,10){$y'$}
\put(20,0){\circle*{3}}
\put(0,20){\circle*{3}}
\put(20,40){\circle*{3}}
\put(40,20){\circle*{3}}
\put(17.7,34){$*$}
\put(2.5,15.5){$*$}
\end{picture}}
\put(95,17){$\longrightarrow$}
\end{picture}
\caption{}\label{figure3}
\end{figure}
is said to be obtained from $\tau$ by
the {\it decorated elementary move} along i.a. $e$.

Denote by $k_t$ and $k^t$, $k=0,1,2$, punctures
and i.a., respectively, on the boundary 
of triangle $t$ with a distinguished corner in accordance with 
Fig.~\ref{figure4}. 
\begin{figure}[h]
\centering
\begin{picture}(60,50)
\put(10,10){\line(1,0){40}}
\put(10,10){\line(1,1){20}}
\put(30,30){\line(1,-1){20}}
\put(10,10){\circle*{3}}
\put(30,30){\circle*{3}}
\put(50,10){\circle*{3}}
\put(43,10){$*$}
\put(28,15){$t$}
\put(28,1){$1^t$}
\put(14,20){$0^t$}
\put(41,20){$2^t$}
\put(28,32){$1_t$}
\put(52,5){$0_t$}
\put(2,5){$2_t$}
\end{picture}
\caption{}\label{figure4}
\end{figure}

With each d.i.t. $\tau$ associate
object
\begin{equation}\label{S-object}
S(\tau):=(\mathbb{R}_+^{\dot{\tau}}\times\mathbb{R}_+^{\dot{\tau}})
((\mathbb{R}_+^{\dot{\tau}}\times\mathbb{R}_+^{\dot{\tau}})/P\mathbb{R}^V_+,
P\mathbb{R}_+^V,\pi_{S(\tau)},\beta_\tau)\in\mathrm{Ob}\mathfrak{P},
\end{equation}
where the free structure group action is defined by
\begin{equation}\label{omega-action}
P\mathbb{R}_+^V\ni f\colon S(\tau)\to
S(\tau),\quad
f^*(\mathbf{t}):=(t_1f(1_t)/f(2_t),t_2f(1_t)/f(0_t));
\end{equation}
\begin{equation}\label{canonical-form}
\beta_\tau:=\sum_{t\in\dot{\tau}}d\ln t_1\wedge d\ln t_2;
\end{equation}
the coordinate functions on $S(\tau)$ used here are defined as follows,
\[
\mathbf{t}:=(t_1,t_2),\quad t_i:=t\circ\mathfrak{pr}_i,\quad t\in\dot{\tau},
\quad \mathfrak{pr}_i\colon S(\tau)\ni (f_1,f_2)
\mapsto f_i\in\mathbb{R}_+^{\dot{\tau}}.
\]
$S(\tau)$ is a Lie group w.r.t the product of
the group structures on $\mathbb{R}_+^{\dot{\tau}}$ factors.

For d.i.t. $\tau$ and $\tau_t$, connected
by the elementary change of decoration in triangle $t$,
see Fig.~\ref{figure2},
 associate mapping
\begin{equation}\label{orientation-map}
S(\tau_t\cdot\tau)\colon S(\tau)\to S(\tau_t),\quad
S(\tau_t\cdot\tau)^*(\mathbf{x}):=\left\{
\begin{array}{cl}
(t_2/t_1,1/t_1)&\mathrm{if}\ x=t';\\
\mathbf{x}&\mathrm{otherwise.}
\end{array}\right.
\end{equation}
For d.i.t. $\tau$ and $\tau^e$, connected
by the decorated elementary move along i.a. $e$ depicted in Fig.~\ref{figure3}, 
associate mapping 
\begin{equation}\label{ptolemy}
S(\tau^e\cdot\tau)\colon S(\tau)\to S(\tau^e),\quad
S(\tau^e\cdot\tau)^*(\mathbf{t}):=\left\{
\begin{array}{cr}
\mathbf{x}\bullet\mathbf{y}& t=x';\\
\mathbf{x}*\mathbf{y}& t=y';\\
\mathbf{t}&\mathrm{otherwise;}
\end{array}\right.
\end{equation}
where 
\begin{equation}\label{dot-star}
\mathbf{x}\bullet\mathbf{y}:=(x_1y_1,x_1y_2+x_2),\quad
\mathbf{x}*\mathbf{y}:=(y_1x_2(x_1y_2+x_2)^{-1},y_2(x_1y_2+x_2)^{-1}).
\end{equation}
Formulae~(\ref{dot-star}) lead to a solution for the pentagon equation, see 
(Kashaev, Sergeev 1996). In fact, they can be obtained by considering
the classical limit of the quantum dilogarithm (Sergeev, 1996).
\begin{proposition}
$S$ extends to a unique covariant functor 
from $\mathfrak{G}_{\tilde{\Delta}_\Sigma}$
into $\mathfrak{P}$.
\end{proposition}
The mapping class group $M_\Sigma$ acts in $\tilde{\Delta}_\Sigma$, 
and for each d.i.t. $\tau$ one has the induced mapping
\[
m\colon\dot{\tau}\to m(\dot{\tau}),\quad \forall m\in M_\Sigma,
\]
of the corresponding set of triangles.
For each $m\in M_\Sigma$ and d.i.t. $\tau$ associate morphism
\[
\mathfrak{s}_m(\tau)\colon S(\tau)\ni f
\mapsto S(\tau\cdot m(\tau))(f\circ m^{-1})\in S(\tau).
\]
\begin{proposition}
Mapping
$ m\mapsto \mathfrak{s}_m$ is a group homomorphism from $M_\Sigma$
into functorial isomorphisms from $S$ to $S$.
\end{proposition}
For a path $\gamma$ in triangle $t\in\dot{\tau}$, 
$\tau\in\tilde{\Delta}_\Sigma$, 
connecting interior points of i.a. $a$ and $b$ with the
initial point being at i.a. $a$, assign function
\begin{equation}\label{homomorphism1}
\mathfrak{u}(\gamma):=\left\{
\begin{array}{cl}
1/t_1& a=1^t,\ b=2^t;\\
t_2 &a=0^t,\ b=1^t;\\
t_2/t_1 &a=0^t,\ b=2^t;
\end{array}\right.\quad \mathfrak{u}(\gamma^{-1}):=1/\mathfrak{u}(\gamma).
\end{equation}
Let $\gamma$ be an element of the first integer homology group of $\Sigma$
represented by an oriented loop $\gamma$ in $\Sigma$ (we shall not
distinguish between loops and the homology classes they define), which
consecutively intersects i.a. $c_1,\ldots,c_{n}$. 
For each $i=1,\ldots,n$ the segment
$\gamma_i$ between i.a.
$c_i$ and $c_{i+1}$ ($c_{n+1}=c_1$) is contained in triangle
$t_i\in\dot{\tau}$. So $\gamma$ can be represented as a composition of paths
\[
\gamma=\gamma_n\cdots\gamma_2\gamma_1,
\]
where definition (\ref{homomorphism1}) for the each factor makes sense.
Consider a sum 
\begin{equation}
\langle\mu_\tau,\gamma\rangle:=\sum_{i=1}^n\ln\mathfrak{u}(\gamma_i)
\end{equation}
which will be called {\it holonomy} aroung $\gamma$.
It is easy to see that this quantity does depend on  
the homology class of $\gamma$ only, and it is linear in the second argument
\[
\langle\mu_\tau,\gamma_1\gamma_2\rangle=
\langle\mu_\tau,\gamma_1\rangle+\langle\mu_\tau,\gamma_2\rangle.
\]
So, it defines {\it momentum mapping}
\begin{equation}\label{momentum-map}
\mu_\tau\colon S(\tau)\to H^1(\Sigma,\mathbb{R}),\quad
\langle\mu_\tau(f),\gamma\rangle:=\langle\mu_\tau,\gamma\rangle(f),
\quad \gamma\in H_1(\Sigma,\mathbb{R}),
\end{equation}
which is in fact a group homomorphism.
\begin{proposition}
The Poisson bracket, associated with symplectic
structure (\ref{canonical-form}),
of holonomies around two loops $\gamma_1$ and $\gamma_2$ is 
given by their intersection index:
$\{\langle\mu_\tau,\gamma_1\rangle,\langle\mu_\tau,\gamma_2\rangle\}
=\gamma_1\circ\gamma_2$.
\end{proposition}
Define a group homomorphism
\begin{equation}\label{hamiltonian}
P\mathbb{R}_+^V\ni f\mapsto 
\xi_f:=\sum_{v\in V}\gamma_v\ln f(v)\in H_1(\Sigma,\mathbb{R}),
\end{equation}
where $\gamma_v$ is a small loop, encircling puncture $v$ in the
counterclockwise direction w.r.t. the orientation of $\Sigma$.
\begin{proposition}
The exponential mapping of the Hamiltonian
vector field, corresponding to function $\langle\mu_\tau,\xi_f\rangle$,
coincides with the group action (\ref{omega-action}).
\end{proposition}
Define a set of group homomorphisms
\begin{equation}\label{psi-map}
\mathfrak{f}(\tau)\colon R(E(\tau))\to S(\tau),\quad
\mathfrak{f}(\tau)^*(\mathbf{t}):=(2^t/1^t,0^t/1^t),\quad\forall\tau\in
\tilde{\Delta}_\Sigma.
\end{equation}
\begin{proposition}
$\mathfrak{f}$ is a functorial morphism from $R\circ E$ to $S$
such that
the following sequence of group homomorphisms is exact 
\[
1\to\mathbb{R}_+\stackrel{\iota}{\longrightarrow}
R\circ E(\tau)\stackrel{\mathfrak{f}(\tau)}{\longrightarrow}
S(\tau) \stackrel{\mu_\tau}{\longrightarrow}
H^1(\Sigma,\mathbb{R})\to0,\quad\forall\tau\in\tilde{\Delta}_\Sigma,
\]
where $\iota(a)(c):=a$, $a\in\mathbb{R}_+$, $c\in E(\tau)$, and
\[
\mathfrak{f}(\tau)\circ\mathfrak{r}_m(E(\tau))=
\mathfrak{s}_m(\tau)\circ\mathfrak{f}(\tau),
\quad\forall m\in M_\Sigma,\quad \forall\tau\in\tilde{\Delta}_\Sigma.
\]
\end{proposition}
To summarize, the Teichm\"uller space $\mathcal{T}_\Sigma$ with the 
Weil--Petersson symplectic structure and the action of the
mapping class group $M_\Sigma$ can be described as the 
Hamiltonian reduction of
$S(\tau)$ w.r.t. $P\mathbb{R}_+^V$ 
over the zero value of the momentum mapping (\ref{momentum-map}):
\begin{equation}\label{hamiltonian-reduction}
\mathcal{T}_\Sigma\simeq\mathbb{R}_+^{E(\tau)}/\mathbb{R}_+^V\simeq
\mu_\tau^{-1}(0)/P\mathbb{R}_+^V,\quad\forall\tau\in\tilde{\Delta}_\Sigma.
\end{equation}

\section{Quantization}\label{section2}

Denote by $\mathfrak{A}$ the category of unital associative algebras
over $\mathbb{C}$. Fix a real positive number $\hbar$.
With each d.i.t. $\tau$ associate object 
$A(\tau)\in\mathrm{Ob}\mathfrak{A}$, which
is the skew-field of fractions of a skew polynomial algebra
generated by elements 
\[
\{\hat{\mathbf{t}}:=(\hat{t}_1,\hat{t}_2)\}_{t\in\dot{\tau}}
\]
subject to the following relations,
\[
\hat{t}_1\hat{t}_2=q^2\hat{t}_2\hat{t}_1,\quad 
\hat{t}_i\hat{t'}_j=\hat{t'}_j\hat{t}_i,\quad i,j=1,2,\ \forall t\ne t'\in
\dot{\tau},\quad q:=\exp(\sqrt{-1}\hbar).
\]
For d.i.t. $\tau$ and $\tau_t$, connected
by the elementary change of decoration in triangle $t$,
see Fig.~\ref{figure2},
 associate morphism
\begin{equation}\label{aorientation-map}
A(\tau_t\cdot\tau)\colon A(\tau_t)\to A(\tau),\quad
A(\tau_t\cdot\tau)(\hat{\mathbf{x}}):=\left\{
\begin{array}{cl}
(q\hat{t}_1^{-1}\hat{t}_2,\hat{t}_1^{-1})&\mathrm{if}\ x=t';\\
\hat{\mathbf{x}}&\mathrm{otherwise;}
\end{array}\right.
\end{equation}
For d.i.t. $\tau$ and $\tau^e$, connected
by the decorated elementary move along i.a. $e$ depicted in Fig.~\ref{figure3}, 
associate morphism
\begin{equation}\label{aptolemy}
A(\tau^e\cdot\tau)\colon A(\tau^e)\to A(\tau),\quad
A(\tau^e\cdot\tau)(\hat{\mathbf{t}})=\left\{
\begin{array}{cr}
\hat{\mathbf{x}}\bullet\hat{\mathbf{y}}& t=x';\\
\hat{\mathbf{x}}*\hat{\mathbf{y}}&t=y';\\
\hat{\mathbf{t}}&\mathrm{otherwise;}
\end{array}\right.
\end{equation}
where notation (\ref{dot-star}) is used.
\begin{proposition}
$A$ extends uniquely to a contravariant functor from 
$\mathfrak{G}_{\tilde{\Delta}_\Sigma}$ into $\mathfrak{A}$.
\end{proposition}
Each $m\in M_\Sigma$ induces algebra isomorphism
\[
m\colon A(\tau)\ni \hat{\mathbf{t}}\mapsto
\hat{\mathbf{t}}'\in A(m(\tau)),\quad t'=m(t)\in m(\dot{\tau}),
\quad\forall\tau\in\tilde{\Delta}_\Sigma.
\]
Define
\[
\mathfrak{a}_m(\tau):=A(m(\tau)\cdot\tau)\circ m\colon A(\tau)\to A(\tau),
\quad\forall m\in M_\Sigma,\quad\forall\tau\in\tilde{\Delta}_\Sigma
\]
\begin{proposition}
Mapping $m\mapsto\mathfrak{a}_m$ is a group homomorphism
from $M_\Sigma$ to functorial automorphisms of $A$.
\end{proposition}
Thus, representations of algebra $A(\tau)$ by linear operators in Hilbert
spaces should lead to (projective) 
representations of the mapping class group $M_\Sigma$.

\subsection{Representations}
\subsubsection{Non-compact representation}
The $*$-algebra structure in $A(\tau)$ defined by
\begin{equation}\label{star-structure}
\hat{t}_i^*=\hat{t}_i,
\end{equation}
is natural from the viewpoint of the underlying classical
phase space. Consider
\[
B(\tau):=
L^2(\mathbb{R}_+^{\dot{\tau}},d\mu(\tau):=\prod_{t\in\dot{\tau}}d\ln t),
\quad \forall\tau\in\tilde{\Delta}_\Sigma,
\]
the Hilbert
space of square integrable functions on $\mathbb{R}_+^{\dot{\tau}}$ w.r.t.
the measure $d\mu(\tau)$. Algebra $A(\tau)$ is realized 
linearly in $B(\tau)$ through the formulae
\[
\mathfrak{b}(\hat{t}_1)f:=ft,\quad
\mathfrak{b}(\hat{t}_2)f:=\exp(-2\sqrt{-1}\hbar t\partial/\partial t)f,\quad
\forall f\in B(\tau).
\]
\begin{proposition}
Let function $\psi(z)$ be a solution of the functional equation
\begin{equation}\label{functional-equation}
\psi(z-\sqrt{-1}\hbar)=\psi(z+\sqrt{-1}\hbar)(1+\exp z),\quad z\in\mathbb{C},
\end{equation}
then operator
\[
T_{x,y}:=y^{-x\partial/\partial x}
\psi(\ln x+2\hbar\sqrt{-1}(x\partial/\partial x-y\partial/\partial y))
\]
satisfies equations
\[
T_{x,y}\mathfrak{b}(\hat{\mathbf{x}}\bullet\hat{\mathbf{y}})=
\mathfrak{b}(\hat{\mathbf{x}})T_{x,y},\quad
T_{x,y}\mathfrak{b}(\hat{\mathbf{x}}*\hat{\mathbf{y}})=
\mathfrak{b}(\hat{\mathbf{y}})T_{x,y}.
\]
\end{proposition}
A particular solution to eqn~(\ref{functional-equation}) is given 
by the formula (Faddeev, 1995)
\[
\psi(z):=\exp
\frac{1}{4}\int_{-\infty}^{+\infty}\frac{\exp(-\sqrt{-1}xz)}{
\sinh(\pi x)\sinh(\hbar x)}\frac{dx}{x},
\]
where the singularity of the integrand at $x=0$ is put below the
contour of integration. In this case operator $T_{x,y}$ is unitary.

\subsubsection{Compact representation}
Let $N\ge2$ be a positive integer, and $\omega$, a complex
primitive $N^{th}$ root of unity.
With any d.i.t. $\tau$ associate object
\[
C(\tau):=\{f\colon S(\tau)\to M(\tau)\}\in \mathrm{Ob}\mathfrak{A},
\]
where unital finite dimensional algebra $M(\tau)$
over $\mathbb{C}$ is generated by elements
\[
\{\check{\mathbf{t}}:=(\check{t}_1,\check{t}_2)\}_{t\in\dot{\tau}}
\]
subject to the following relations,
\[
\check{t}_1\check{t}_2=\omega\check{t}_2\check{t}_1,\quad 
\check{t}_i\check{t'}_j=\check{t'}_j\check{t}_i,\quad 
\check{t}_i^N=1,\quad i,j=1,2,\ \forall t\ne t'\in \dot{\tau};
\]
the algebra structure in $C(\tau)$ being given by the pointwise algebraic
operations in $M(\tau)$. Note that algebra $M(\tau)$ has a unique
up to isomorphism irreducible finite dimensional representation
with the $*$-algebra structure given by
\[
\check{t}_i^*=\check{t}_i^{-1},
\quad i=1,2,\quad\forall t\in\dot{\tau}.
\]
Algebra $A(\tau)$ at $q^2=\omega$ is represented in $C(\tau)$ through
the formulae
\[
\mathfrak{c}(\tau)\colon
A(\tau)\to C(\tau),\quad
\mathfrak{c}(\tau)(\hat{t}_i)(h):=\check{t}_i\sqrt[N]{t_i(h)},\quad
i=1,2,\ \forall h\in S(\tau),
\]
where positive $N^{th}$ roots are taken. This repsentation does not respect
the $*$-structure (\ref{star-structure}).

For d.i.t. $\tau$ and $\tau'$, where $\tau'=\tau_t$ (elementary change
of decoration in triangle $t$) or $\tau'=\tau^e$ (elementary move along
i.a. $e$), associate morphisms
\[
C(\tau'\cdot\tau)\colon C(\tau')\to C(\tau),
\]
\[
C(\tau'\cdot\tau)(f)(h):=M_h(\tau'\cdot\tau)(f(S(\tau'\cdot\tau)(h))),\quad
\forall f\in C(\tau'),\quad \forall h\in S(\tau),
\]
where algebra isomorphisms 
\[
M_h(\tau'\cdot\tau)\colon M(\tau')\to M(\tau),\quad \forall h\in S(\tau)
\]
are defined by
\[
M_h(\tau_t\cdot\tau)(\check{\mathbf{x}}):=\left\{
\begin{array}{cl}
(\omega^{1/2}\check{t}_1^{-1}\check{t}_2,\check{t}_1^{-1})&
\mathrm{if}\ x=t';\\
\check{\mathbf{x}}&\mathrm{otherwise,}
\end{array}\right.\quad \omega^{N/2}=(-1)^{N-1};
\]
\[
M_h(\tau^e\cdot\tau)(\check{\mathbf{t}}):=\left\{
\begin{array}{cr}
\check{\mathbf{x}}\bullet_h\check{\mathbf{y}}& t=x';\\
\check{\mathbf{x}}*_h\check{\mathbf{y}}&t=y';\\
\check{\mathbf{t}}&\mathrm{otherwise;}
\end{array}\right.
\]
\[
\check{\mathbf{x}}\bullet_h\check{\mathbf{y}}:=
\left(\check{x}_1\check{y}_1,
(\check{x}_2+\check{x}_1\check{y}_2h_{x,y})
/\sqrt[N]{1+h_{x,y}^N}\right),\quad h_{x,y}:=\sqrt[N]{x_1y_2/x_2}(h),
\]
\[
\check{\mathbf{x}}*_h\check{\mathbf{y}}:=
(\check{y}_1\check{x}_2(\check{\mathbf{x}}\bullet_h\check{\mathbf{y}})_2^{-1},
\check{y}_2(\check{\mathbf{x}}\bullet_h\check{\mathbf{y}})_2^{-1}),
\]
with positive $N^{th}$ roots being chosen;
see also Figs~\ref{figure2},~\ref{figure3} for the notation of the triangles.

\begin{proposition}
$C$ extends uniquely to a contravariant functor from 
$\mathfrak{G}_{\tilde{\Delta}_\Sigma}$ into $\mathfrak{A}$,
$\mathfrak{c}$ being a functorial morphism from $A$ to $C$.
\end{proposition}
\begin{proposition}
Let $\Psi_{\omega,\lambda}(w)$ be a solution of the functional
equation
\begin{equation}\label{functional-equation-compact}
\Psi_{\lambda}(\omega w)(1-w\lambda)/\lambda'=
\Psi_{\lambda}(w),
\end{equation}
where
\[
\lambda,\lambda'\in\mathbb{R}_+,\ 
\lambda'=\sqrt[N]{1+\lambda^N},\quad w^{N}=-1;
\]
then the element
\[
T_{h,x,y}:=(\sum_{i,j=0}^{N-1}\omega^{-ij}\check{y}_1^i\check{x}_2^j)
\Psi_{h_{x,y}}(-\check{x}_2^{-1}\check{x}_1\check{y}_2),
\]
satisfies the equations
\[
T_{h,x,y}\check{\mathbf{x}}\bullet_h\check{\mathbf{y}}=
\check{\mathbf{x}}T_{h,x,y},\quad
T_{h,x,y}\check{\mathbf{x}}*_h\check{\mathbf{y}}=
\check{\mathbf{y}}T_{h,x,y}.
\]
\end{proposition}
Solution of the functional equation (\ref{functional-equation-compact})
is unique up to a complex factor.
It satisfies the pentagon equation (Faddeev, Kashaev 1992)
which is a non-commutative analogue of Roger's five-term identity 
for the dilogarithm.

\section*{Summary}

The Teichm\"uller space of 
marked conformal types of hyperbolic metrics
on a punctured surface with 
the Weil--Petersson symplectic structure
and the action of the mapping class group can be described as the Hamiltonian
reduction (\ref{hamiltonian-reduction}) of a finite dimensional symplectic
manifold (\ref{S-object}). The quantization of the latter is straightforward,
and the action of the mapping class group is realized through the quantum
dilogarithms. These results relate the quantum hyperbolic invariant of
knots (Kashaev 1994, 1995, 1997) 
to the quantum theory of Teichm\"uller spaces of punctured surfaces.

\end{document}